\documentclass{kjmath}

\usepackage{amssymb}
\usepackage{eucal}
\usepackage{graphicx}
\usepackage{color}
\usepackage{url}
 \usepackage{tikz-cd}

\begin{document}

\kjmvolume{%
@VOL@
}
\kjmnumber{%
@NUM@
}
\kjmyear{%
@YEAR@
}
\kjmstpage{%
@STPAGE@
}
\kjmendpage{%
@ENDPAGE@
}


\title[Functor Of Relations on Hypergroups and Hypermodules]{Functors Associated to Relations on Hypergroups and Hypermodules}


\author[B.\ Afshar]{Behnam Afshar$^1$}

\address{$^1$Department of Mathematics, Statistics and Computer Science,
	\newline \indent Faculty of Science,
	\newline \indent University of Tehran, Tehran, Iran}
\email{behnamafshar@ut.ac.ir
\newline \indent ORCID iD: \url{https://orcid.org/0009-0006-7095-5795}}

\author[R.\ Ameri]{Reza Ameri$^2$}


\address{$^2$Department of Mathematics, Statistics and Computer Science,
\newline \indent Faculty of Science,
\newline \indent University of Tehran, Tehran, Iran}
\email{rameri@ut.ac.ir
\newline \indent ORCID iD: \url{https://orcid.org/0000-0001-5760-1788}}



\begin{abstract}
If $H$ is a strongly regular hypergroup, we show that the set of regular relations on $H$ and the set of subhypergroups containing $0_{H}$ are two lattices that are isomorphic to each other. In the next step, we introduce and study the properties of functors that are constructed by a sequence of strongly regular relations. This helps us to define a specific type of free objects and tensor products on the category of regular hypergroups.

\end{abstract}

\maketitle

\section{Introduction} \label{sec:introduction}

The fundamental relations are one of the most important and interesting concepts in algebraic hyperstructures that ordinary algebraic structures are derived from algebraic hyperstructures by them. The fundamental relation $\beta^*$ on hypergroups was defined by M. Koskas, P. Corsini, D. Ferni and T. Vogiouklis.
Then D. Ferni introduced the fundamental relation $\gamma^{*}$ which is the transitive closure of $\gamma$ and is the smallest relation such that $H/\gamma^*$ is an abelian group.
Subsequently, fundamental relations were gradually introduced on other algebraic hyperstructures. In \cite{Afshar-Ameri 1}, it has
been demonstrated that relations $\beta^*$ and $\gamma^{*}$ are related together in the form of $\gamma^{*}=\delta\ast\beta^*$, where $\delta$ is the congruence relation with respect to the commutator subgroup. In \cite{Afshar-Ameri 2}, a one-to-one correspondence is also introduced between strongly regular relations on regular hypergroups and subhypergroups containing $S_{\beta}$. Then, this one-to-one correspondence is elevated to a lattice isomorphism.\\
In Section \ref{Sec3}, we will prove this correspondence and isomorphism for regular relations on canonical hypergroups. In Section \ref{Sec4}, using sequences of strongly regular relations, we introduce a new type of functor on the category of hypergroups, which we call the \emph{functor of relations}. Then in Section \ref{Sec5}, we will give a modular structure to the category of hypergroups and introduce and study a type of tensor product and the concept of free objects on it using the functor of relations. One useful application of the functor of relations is that, based on a specific property, one can construct the corresponding functor by making appropriate choices of relations.

\section{Preliminaries} \label{subsec:SecUnit}

\begin{definition}[\cite{Corsini}]
	A hypergroup $H$ is called a regular hypergroup if $E_{H}\neq\emptyset$ and $C_{H}(x)\neq\emptyset$ for every $x\in H$. A regular hypregroup $H$ is called a strongly regular hypergroup, if $C_{H}(x)=\{x^{-1}\}$ for every $x\in H$.
\end{definition}

If $H$ is a strongly regular hypergroup, then $|E_{H}|=1$. Because for every $e_{1},e_{2}\in E_{H}$; $e_{1}\in e_{2}\circ e_{1}$ and $e_{2}\in e_{2}\circ e_{2}$. So $e_{1},e_{2}\in C_{H}(e_{2})$ and $e_{1}=e_{2}$. A subhypergroup $K$ of $H$, is closed if and only if $0_{H}\in K$.

\begin{corollary}[\cite{Davvaz}]
	If $(H,\circ)$ is a hypergroup and $\mathcal{R}$ be an equivalence relation on $H$, then $\mathcal{R}$ is regular (strongly regular) if and only if $(H/\mathcal{R},\otimes)$ is a hypergroup (group).
\end{corollary}

\begin{definition}[\cite{Davvaz}]
	For all $n>1$, we define the relations $\beta_{n}$ and $\gamma_n$ on a semi-hypergroup $H$, as follows:
	\begin{center}
		$a\beta_{n}b\Longleftrightarrow \exists(x_{1},x_{2},...,x_{n})\in H^{n} : \{a,b\}\subseteq\prod_{i=1}^{n}x_{i}$\ ;
	\end{center}
	\begin{center}
		$ a\gamma_{n}b\Longleftrightarrow \exists(x_{1},x_{2},...,x_{n})\in H^{n} , \sigma\in \mathbb{S}_{n}: \ a\in\prod_{i=1}^{n}x_{i}  \ , \  b\in\prod_{i=1}^{n}y_{\sigma(i)} $,
	\end{center}	
	and $\beta=\bigcup_{n\geq1}\beta_{n}$ and $ \gamma=\bigcup_{n\geq1}\gamma_{n}$, where $\beta_{1}=\gamma_{1}=\{(x,x);x\tiny{\in} H\}$.
	Also $\beta^{*}$ and $\gamma^{*}$ are transitive closures of $\beta$ and $\gamma$ respectively.
\end{definition}

If $H$ is a hypergroup then $\gamma=\gamma^{*}$ and $\beta=\beta^{*}$. Also $\beta^{*}$ is the smallest strongly regular relation on $H$ and $\gamma^{*}$ is the smallest strongly regular relation on $H$  such that the quotient $H/\gamma^{*}$ is commutative group \cite{Corsini}. Also, consider that $\eta_{a}^{\ast}$ and $\tau_{n}^{\ast}$ as introduced in \cite{Davvaz}, are the smallest strongly regular relations such that $H/\eta_{a}^{\ast}$ and $H/\tau_{n}^{\ast}$ are, respectively, a cyclic group and a solvable group.

Let $SR(H)$ be the set of all strongly regular relations on strongly regular hypergroup $H$, and $S(H):=\{S_{\rho}\ ;\  \rho\in SR(H)\}$ where $S_{\rho}:=\{s\in H \ ;\ \rho(s)=e_{H/\rho}\}$. Also, denote the set of all subhypergroups of $H$ containing $S_{\beta}$ by $(S_{\beta})$, and the set of all normal subhypergroups of $H$ containing $S_{\beta}$ by $N(S_{\beta})$. In \cite{Afshar-Ameri 2}, it is proved that $S(H)\cong N(S_{\beta})$ is a lattice isomorphism.

\begin{proposition}[\cite{Jantosciak}]\label{Jantosciak}
	In a transposition hypergroup, if $\rho$ is a regular equivalence
	relation that has a nonempty set $N$ of scalar identities, then $N$ is a
	reflexive closed set and $\rho$ is the relation equivalence modulo $N$.
\end{proposition}

A hypergroup is a polygroup if and only if it is a transposition
hypergroup with a scalar identity, and a join space is a commutative transposition hypergroup, \cite{Jantosciak}. A subhypergroup $N$ of hypergroup $H$ is normal if $aN=Na$, for all $a\in H$.

\begin{proposition}[\cite{Jantosciak}]\label{Jantosciak2}
	In a transposition hypergroup, a normal closed set is
	reflexive.
\end{proposition}

Let $\mathcal{A}$ be a subcategory of a category $\mathcal{B}$, then $\mathcal{A}$ is a \emph{full subcategory} of $\mathcal{B}$ if for each pair of objects $X$ and $Y$ of $A$; $Hom_{\mathcal{A}}(X,Y)=Hom_{\mathcal{B}}(X,Y)$. A functor $F$ from category $\mathcal{C}$ to category $\mathcal{D}$ is said to be \emph{full} if $F_{X,Y}$ is surjective, and is said to be \emph{faithful} if $F_{X,Y}$ is injective. Where $F_{X,Y}$ is induced function of $F$,
\begin{center}
	$F_{X,Y}:Hom_{\mathcal{C}}(X,Y)\rightarrow Hom_{\mathcal{D}}(F(X),F(Y))$.
\end{center}
on hom-sets, \cite{Awodey}.

\section{Regular relations and closed subhypergroups}\label{Sec3}

In this section $(H,+)$ is a canonical hypergroup. We denote the set of all regular relations on $H$ by $R(H)$ and consider $S_{H}=\{S_{\rho}\ ;\ \rho\in R(H)\}$, where:
\begin{equation}\label{eq1}
	S_{\rho}:=\{s\in H \ ; \ \rho(x)\in \rho(x)\oplus\rho(s), \ \forall x\in H \}
\end{equation}
and $\rho(x)\oplus\rho(s)=\{\rho(y);\ y\in x+s\}$.
A subhypergroup of $H$, is closed if and only if it contains $0_{H}$  \cite{Massouros}. Since $H/\rho$ is a canonical hypergroup, then  $0_{H}\in S_{\rho}$. Also, we denote the set of all closed subhypergroups of $H$ (equivalently, subhypergroups of $H$ containing $0_{H}$), by $(0_{H})$.

\begin{lemma}        
	If $H$ is a canonical hypergroup, then $S_{H}\subseteq (0_{H})$.
	\begin{proof}
		Let $x\in S_{\rho}$. Since $S_{\rho}=\rho(0_{H})$, then:
		\begin{equation*}	
			\rho(x+\rho(0))=\rho(x)\oplus\rho(\rho(0))=\rho(x)\oplus\rho(0)=\rho(0).
		\end{equation*}	
		Hence $x+S_{\rho}\subseteq S_{\rho}$. Conversely, suppose that $s\in \rho(0)$. Since $H=x+H$, then there exists $h\in H$ such that $s\in x+h$ and $\rho(s)\in \rho(x)\oplus\rho(h)$. From that $x,s\in S_{\rho}=\rho(0)$, then $\rho(0)\in \rho(0)\oplus\rho(h)=\rho(h)$ and $h\in\rho(0)$. Therefore $s\in x+S_{\rho}$ and $S_{\rho}\subseteq x+\rho(0)$.

	\end{proof}	
\end{lemma}

Therefore, $S_{\rho}$ is a closed (canonical) subhypergroup of $H$. Because $S_{\rho}$ is normal, then by Proposition \ref{Jantosciak2}, it is reflexive.

\begin{lemma}         \label{2.1}
	If $\rho\in R(H)$, then $\rho(x)=x+S_{\rho}$ for every $x\in H$.
	\begin{proof}
		This is a straightforward results of Proposition \ref{Jantosciak}.	
	\end{proof}	
\end{lemma}

Clearly, $S_{\rho}=\rho(0_{H})\in (0_{H})$, and the congruence relation
modulo of each element of $(0_{H})$ is regular. So $S_{H}=(0_{H})$.

\begin{theorem}             \label{2.2}
	There is a bijection between the sets $R(H)$ and $(0_{H})$.
	\begin{proof}   
		It is enough to consider the maps $f$ and $g$ as follows:
		\begin{equation}\label{eq2}
			\underset{\theta\ \longmapsto \ S_{\theta}}{f:R(H)\rightarrow (0_{H})} \ \ \ \ \  and \ \ \ \ \  \underset{S \ \longmapsto \ \theta_{S}}{g:(0_{H})\rightarrow R(H)}
		\end{equation}
		where $\theta_{S}=\{(x,y)\ ;\ x+S=y+S\}$. The congruence relation modulo of each element of $(0_{H})$ is regular. So $g$ is well defined. If $S\in (0_{H})$, then $\theta_{S}(0)=S$ and $H/\theta_{S}=H/S$ is an isomorphism of canonical hypergroups. So:
		\begin{align*}
			f\circ g (S)&=f\big{(}\{(x,y);\ x+S=y+S\}\big{)}\\
			&=\{h\in H;\ \theta_{S}(x)\oplus\theta_{S}(h)=\theta_{S}(x),\ \forall x\in H\}\\
			&=\{h\in H;\  \theta_{S}(h)=\theta_{S}(0)\}\\
			&=\{h\in H;\  h+S=S\}\\
			&=S.
		\end{align*} 
		Conversely, if $\theta\in R(H)$, then:
		\begin{center}
			$g\circ f(\theta)=g(S_{\theta})=\{(x,y;\ x+S_{\theta}=y+S_{\theta})\}=\{(x,y);\ \theta(x)=\theta(y)\}=\theta$.
		\end{center}
	\end{proof}
\end{theorem}

Therefore  $H/\rho\cong H/S_{\rho}$ is a good isomorphism of hypergroups, for every $\rho\in R(H)$.
For a hypergroup $H$, $(S_{H},\subseteq)$ and $(R(H),\subseteq)$ are posets and $S_{H}=\{S_{\theta};\ \theta\tiny{\in}R(H)\}=(0_{H})$, $R(H)=\{\theta_{S};\ S\tiny{\in}(0_{H})\}$.\\
Let $\Delta_{H}=\{(h,h); h\tiny{\in}H\}$ and $\nabla_{H}=H^{2}$. Then $\nabla_{H},\Delta_{H}\in R(H)$. For every $\rho,\sigma\in R(H)$ put $\rho\tiny{\vee}\sigma=\rho\cup(\rho \tiny{\circ}\sigma)\cup(\rho\tiny{\circ}\sigma \tiny{\circ}\rho)\cup(\rho\tiny{\circ}\sigma\tiny{\circ}\rho\tiny{\circ}\sigma)\cup...$ where $\rho \tiny{\circ}\sigma=\{(x,y)\tiny{\in} H^{2};\ \exists z\tiny{\in} H \ni (x,z)\tiny{\in}\rho, (z,y)\tiny{\in}\sigma\}$. Also $\rho \cap\sigma\in R(H)$ and $\bigcap_{\rho\in R(H)}\rho=\Delta_{H}$.\\
If $L$ is a lattice and $a\in L$ then $I(a):=\{b\tiny{\in}L;\ b\leqslant a\}$ is principal ideal generated by $a$ and $F(a):=\{b\tiny{\in}L;\ a\leqslant b\}$ is filter generaed by $a$, \cite{Burris}.

\begin{lemma}                  \label{3.1}
	If $\rho,\sigma\in R(H)$, then $\rho\tiny{\vee}\sigma\in R(H)$.
	\begin{proof}
		Clearly, $\rho\tiny{\vee}\sigma$ is the smallest equivalence relation containing $\rho$ and $\sigma$. Let $(x,y)\in \rho\tiny{\vee}\sigma$ and $z\in H$. Then there are $n\in \mathbb{N}$ and $z_{1},...,z_{n}\in H$ such that $z_{1}=x$, $z_{n}=y$ and $(z_{i},z_{i+1})\in \rho\cup\sigma$, for every $1\leq i\leq n-1$. Without loss of generality assume that $(x,z_{2})\in \rho$, $(z_{2},z_{3})\in \sigma$,..., $(z_{n-1},y)\in \rho$. So
		\begin{center}
			$z+x\ \bar{\rho}\ z+z_{2}\ \bar{\sigma}\ z+z_{3}\ ...\ z+z_{n-1}\ \bar{\rho}\ z+y.$
		\end{center}
		Hence for every $t_{1}\in z+x$ there are $t_{2},...,t_{n}\in H$, where $t_{k}\in z+z_{k}$ ($2\leq k\leq n$), and $t_{1}\ \rho\ t_{2}\ \sigma\ t_{3}\ ... \ t_{n-1}\ \rho\ t_{n}$.
		Thus $(t_{1},t_{n})\in\rho\vee\sigma $ and $z+x\ \overline{\rho\tiny{\vee}\sigma}\ z+y$.	 	
	\end{proof}	
\end{lemma}

Therefore, it is easy to see that $S_{H}$ and $R(H)$ are complete lattices, and
\begin{center}
	$\bigvee_{\rho\in R(H)}\rho=\nabla_{H}$,\ \  $\bigcap_{\rho\in R(H)}\rho=\Delta_{H}$,\ \  $\bigvee_{S\in S_{H}}S=H$,\ \   $\bigcap_{S\in S_{H}}S=0_{H}$.
\end{center}

\begin{proposition} \label{3.3}
	The lattices $S_{H}$ and $R(H)$ are isomorph.
	\begin{proof}
		If $\rho\subseteq\theta$, then $\rho(0)\subseteq\theta(0)$. So the function $f$ in \eqref{eq2} is order preserving. 
	\end{proof}		
\end{proposition}

Therefore $S_{\rho}+S_{\theta}=S_{\rho\vee\theta}$, $S_{\rho}\cap S_{\theta}=S_{\rho\cap\theta}$, $\theta_{S}\vee \theta_{T}=\theta_{S+T}$ and $\theta_{S}\cap \theta_{T}=\theta_{S\cap T}$. Where $S,T\in S_{H}$, $\rho,\sigma\in R(H)$ and $\theta_{S}=\{(x,y)\in H^{2};\ x+S= y+S\}$. Also $\bigvee_{\rho\in R(H)}\rho=\nabla_{H}$, $\bigcap_{\rho\in R(H)}\rho=\Delta_{H}$, $\bigvee_{S\in S_{H}}S=H$ and $\bigcap_{S\in S_{H}}S=0_{H}$.

\begin{theorem} \label{T3.6}
	Let $f:H\rightarrow K$ be a good homomorphism of canonical hypergroups and $\rho\in R(H)$, $\sigma\in R(K)$ such that $f(S_{\rho})\leqslant S_{\sigma}$. Then there is a unique good homomorphism $\bar{f}:H/\rho\rightarrow K/\sigma$	where $\bar{f}(\rho(h))=\sigma(f(h))$. Moreover, $\bar{f}$ is a isomorphism if and only if $Imf\vee S_{\rho}=K$ and $f^{-1}(S_{\sigma})\subseteq S_{\rho}$.
	\begin{proof}
		If $\rho(x)=\rho(y)$, then $(f(x),f(y))\in f(S_{\rho})\subseteq S_{\sigma}$ and $\sigma(f(x))=\sigma(f(y))$. Thus $\bar{f}$ is welldefined. Also $\bar{f}(\rho(x)\oplus\rho(y))=\bar{f}(\rho(x+y))=\sigma(f(x+y))=\sigma(f(x)+f(y))=\sigma(f(x))\oplus\sigma(f(y))=\bar{f}(\rho(x))\oplus\bar{f}(\rho(y))$. Let $\bar{f}$ be an isomorphism, then $f$ is onto and $Imf\vee S_{\sigma}=K$. Since $Ker\bar{f}=S_{\rho}$ then $\{\rho(x);\ \bar{f}(\rho(x))=S_{\sigma}\}=\{\rho(x);\ \sigma(f(x))=S_{\sigma}\}=\{x+S_{\rho};\ f(x)\tiny{\in} S_{\sigma}\}=S_{\rho}$. Thus, $f^{-1}(S_{\sigma})\subseteq S_{\rho}$. If $f^{-1}(S_{\sigma})\subseteq S_{\rho}$ and $Imf\vee S_{\sigma}=K$, then $Ker\bar{f}=S_{\rho}$ and $S_{\sigma}\subseteq Imf$. Therefore, $Imf=K$ and $\bar{f}$ is onto. 
	\end{proof}
\end{theorem}

\begin{example} \label{3--7}
	Consider the canonical hypergroup $H$ as follows:
	\begin{center}
		\begin{tabular}{c|cccccccc}		
			$+$& $0$  & $x$      &$y$  &$z$  &$u$  &$v$ &$-x$ & $-z$\\
			\hline
			$0$& $0$  & $x$        &$y$       &$z$      &$u$      &$v$    &$-x$     & $-z$\\
			$x$& $x$  & $y$,$v$    &$z$,$-x$   &$u$      &$x$,$-z$  &$-x$    &$0$,$u$ & $y$\\
			$y$& $y$  & $z$,$-x$    &$0$,$u$   &$x$      &$y$,$v$  &$u$    &$x$,$-z$ & $-x$\\
			$z$& $z$  & $u$        &$x$       &$v$      &$-x$      &$-z$    &$y$     & $0$\\
			$u$& $u$  & $x$,$-z$    &$y$,$v$   &$-x$      &$0$,$u$  &$y$    &$z$,$-x$ & $x$\\
			$v$& $v$  & $-x$        &$u$       &$-z$      &$y$      &$0$    &$x$     & $z$\\
			$-x$& $-x$  & $0$,$u$    &$x$,$-z$   &$y$      &$z$,$-x$  &$x$    &$y$,$v$ & $u$\\
			$-z$& $-z$  & $y$        &$-x$       &$0$      &$x$      &$z$    &$u$     & $v$
		\end{tabular}
	\end{center}
	Then $\gamma=\beta=\{(y,v),(v,y),(z,-x),(-x,z),(x,-z),(-z,x),(0,u),(u,0)\}\cup\Delta_{H}$, and $S_{\gamma}=S_{\beta}=\{0,u\}$.	Also, $\rho=\Delta_{H}\cup S^{2}\cup T^{2}$, where $S=\{0,u,v,y\}$ and $T=\{x,-x,z,-z\}$.
	\begin{figure}[htp]
		\centering
		\includegraphics[scale=0.4]{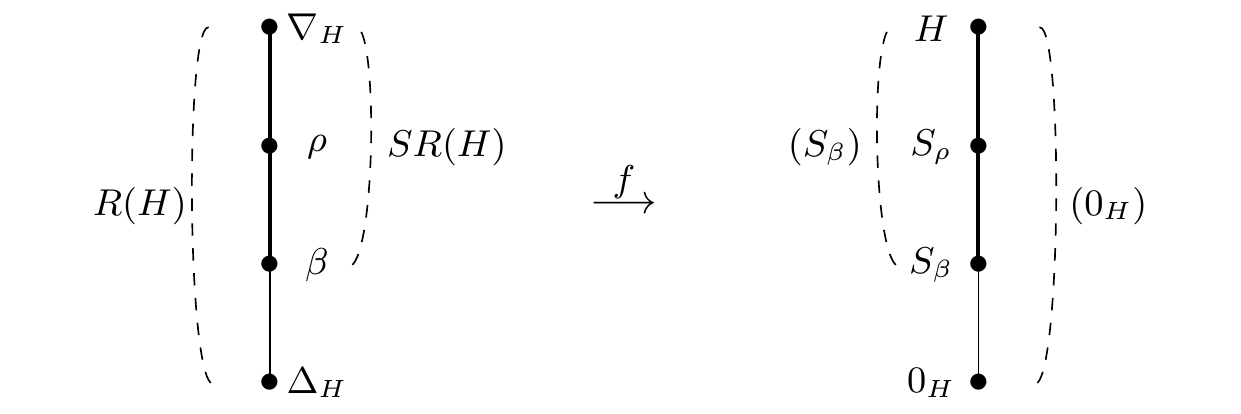}
	\end{figure}
	Therefore $S_{\rho}=\{0,u,v,y\}$, and $f|_{SR(H)}$ is an isomorphism of lattices of strongly regular relations on $H$ and subhypergroups of $H$ containing $S_{\beta}$. Also $H/S_{\beta}\cong \mathbb{Z}_{4}$ and $H/S_{\rho}\cong \mathbb{Z}_{2}$.
\end{example}

If $\rho\in R(H)$, then by Theorem \ref{3.3} and results in \cite{Afshar-Ameri 2}, we have:
\begin{align*}
	\rho\in SR(H)&\Leftrightarrow S_{\rho} \ is \ complete\\
	&\Leftrightarrow S_{\beta} \subseteq S_{\rho}\\
	&\Leftrightarrow \beta\subseteq \rho.
\end{align*}
Therefore $SR(H)=\beta\vee R(H)=\{\beta\vee\rho \ ;\ \rho\in R(H)\}$.

\section{Functor of relations}\label{Sec4}

Let $f:H\rightarrow K$ be a good homomorphism of strongly regular hypergroups and $\rho\in SR(H)$, $\sigma\in SR(K)$ such that $f(S_{\rho})\leq S_{\sigma}$. Then there is
a unique homomorphism $\bar{f}:H/\rho \rightarrow K/\sigma$ where $\bar{f}(\rho(h))=\sigma(f(h))$, \cite[Corollary 4.5]{Afshar-Ameri 2}.

\begin{proposition}\label{P1}
	If $f:H\rightarrow K$ is an inclusion homomorphism of regular hypergroups, and $\rho\in SR(H)$, $\sigma\in SR(K)$, such that $f(S_{\rho})\leq S_{\sigma}$, then there is
	a unique homomorphism of groups $\bar{f}:H/\rho \rightarrow K/\sigma$, where $\bar{f}(\rho(h))=\sigma(f(h))$.
	\begin{proof}	
		Let $\rho(x)=\rho(y)$, then by \cite[Lemma 3.5]{Afshar-Ameri 2},  $x+S_{\rho}=y+S_{\rho}$ and $f(x+S_{\rho})=f(y+S_{\rho})$. So $f(x+S_{\rho})\subseteq f(x)+f(S_{\rho})\subseteq f(x)+S_{\sigma}=\sigma(f(x))$ and $f(x+S_{\rho})\subseteq f(y)+f(S_{\rho})\subseteq f(y)+S_{\sigma}=\sigma(f(y))$. Therefore $\sigma(f(x))=\sigma(f(y))$, and $\bar{f}$ is well defined. Also
		$\bar{f}(\rho(x)\oplus\rho(y))=\bar{f}(\rho(x+y))=\sigma(f(x+y))\subseteq\sigma(f(x)+f(y))=\sigma(f(x))\oplus\sigma(f(y))=\bar{f}(\rho(x))\oplus\bar{f}(\rho(y))$.	
		
	\end{proof}
	
\end{proposition}

Using Lemma \ref{2.1}, Proposition \ref{P1} can be stated for regular relations on canonical hypergroups. That is, if $f:H\rightarrow K$ is an inclusion homomorphism of canonical hypergroups and $\rho\in R(H)$, $\sigma\in R(K)$, such that $f(S_{\rho})\leq S_{\sigma}$, then there is
a unique homomorphism of groups $\bar{f}:H/\rho \rightarrow K/\sigma$, where $\bar{f}(\rho(h))=\sigma(f(h))$.\\

\begin{definition}
	Let $\textbf{H}$ and $\overline{\textbf{H}}$, are categories of regular hypergroups, which their morphisms are good and inclusion homomorphisms, respectively. A sequence $\boldsymbol\theta=\big\{ \theta_{_{H}}\big\}_{H\in I}$ of relations, where $I=Ob(\textbf{H})$ and $\theta_{_{H}}\in R(H)$ for every $H\in I$, is called \emph{hom-invariant} if for every hypergroups $H$, $K$ and $f\in Hom_{\textbf{H}}(H,K)$, (or $f\in Hom_{\overline{\textbf{H}}}(H,K)$), if $(h_{1},h_{2})\in \theta_{H}$, then $(f(h_{1}),f(h_{2}))\in \theta_{K}$.
\end{definition}

Where it dose not lead to ambiguity, we use $S_{H}$ instead of $S_{\theta_{H}}$.
It can be easily investigated that, $\boldsymbol\theta$ is hom-invariant if and only if $f(S_{H})\subseteq S_{K}$, for every homomorphism $f:H\rightarrow K$ of hypergroups.

\begin{theorem} \label{T.1}
	Every hom-invariant sequence $\boldsymbol\theta$ of strongly regular relations, is a functor from $\textbf{H}$ (or $\overline{\textbf{H}}$) to the category of groups $\textbf{Group}$.
	
	\begin{proof}	
		Let $\boldsymbol\theta(H)=H/\theta_{_{H}}$, and for every inclusion (resp. good) homomorphism $f:H\rightarrow K$:
		\begin{equation}
			\underset{\theta_{_{H}}(h)\ \mapsto \ \theta_{_{K}}(f(h)) }{\boldsymbol\theta f:H/\theta_{_{H}}\rightarrow K/\theta_{_{K}}}
		\end{equation}
		Since $\boldsymbol\theta$ is hom-invariant, then $\boldsymbol\theta f\in Hom(H/\theta_{_{H}},K/\theta_{_{K}})$ and $\boldsymbol\theta g\in Hom(K/\theta_{_{K}},T/\theta_{_{T}})$, where $g:K\rightarrow T$ is an inclusion (resp. good) homomorphism of hypergroups. Also $(g\circ f)(S_{H})=g(f(S_{H}))\subseteq S_{T}$ and $\boldsymbol\theta (g\circ f), \boldsymbol\theta g\circ \boldsymbol\theta f \in Hom(H/\theta_{_{H}},T/\theta_{_{T}})$. Moreover, $\boldsymbol\theta (g\circ f)(\theta_{_{H}}(h_{1}))=\theta_{_{T}}((g\circ f)(h_{1}))=\boldsymbol\theta g(\theta_{_{K}}(f(h_{1})))=(\boldsymbol\theta g\circ \boldsymbol\theta f)(\theta_{_{H}}(h_{1}))$.
	\end{proof}	
\end{theorem}

\begin{example} \label{EX6}
	For any hypergroup $H$, consider:
	\begin{equation}
		(a,b)\in\alpha_{H}\Longleftrightarrow  \exists n\in N, (x_{1},x_{2},...,x_{n})\in H^{n}, \sigma\in \mathbb{A}_{n}: \ a\in\prod_{i=1}^{n}x_{i}  \ , \  b\in\prod_{i=1}^{n}y_{\sigma(i)}
	\end{equation}
	and let $\alpha_{H}^{*}$ to be the transitive closure of $\alpha_{H}$. If $\boldsymbol\alpha=\{\alpha_{H}^{*}\}_{H\in I}$, then for every $f\in Hom(H,K)$, $(x,y)\in \alpha_{H}^{*}$ results that $(f(x),f(y))\in \alpha_{K}^{*}$. So $f(S_{H})\subseteq S_{K}$. Therefore $\boldsymbol\alpha$ is hom-invariant sequence. Also $\boldsymbol\beta\subseteq\boldsymbol\alpha\subseteq\boldsymbol\gamma$.	
\end{example}

\begin{example}	\label{EX5}
	Let $\mathfrak{R}=\mathbb{Z}\times \mathbb{Z}$ and $(a,b)\oplus(c,d)=\big\{(a+c,x);\ x\in\mathbb{Z} \big\}$. Then $\mathfrak{R}$ is a regular hypergroup, that $E_{\mathfrak{R}}=S_{\beta}=S_{\gamma}=\{(0,x);\ x\in \mathbb{Z}\}$ and $C\big((a,b)\big)=\{(-a,x);\ x\in \mathbb{Z}\}$. Let $\boldsymbol\xi=\{\xi_{H}\}_{H\in I}$, such that $\xi_{H}=\gamma_{H}$ ($H\neq\mathfrak{R}$), and for a fixed $n\geq2$ consider the relation $\xi_{\mathfrak{R}}$ on $\mathfrak{R}$ as follows:
	\begin{equation}\label{eq4}
		(a,b)\xi_{\mathfrak{R}}(c,d)\Longleftrightarrow a+n\mathbb{Z}=c+n\mathbb{Z}
	\end{equation}
	Clearly $\xi_{\mathfrak{R}}\in SR(\mathfrak{R})$, $\mathfrak{R}/\xi_{\mathfrak{R}}\cong \mathbb{Z}_{n}$ and  $\boldsymbol\xi$ is a sequence of strongly regular relations on hypergroups. But $\boldsymbol\xi$ is not hom-invariant, because for $f\in Hom(\mathfrak{R},\mathbb{Z})$ such that $f((r,s))=r$, we have $f(\{(r,s); r\in n\mathbb{Z}\})=n\mathbb{Z}\nsubseteq\{0\}$, and
	\begin{equation*}
		\underset{\xi_{\mathfrak{R}}(r,s)\mapsto r}	
		{\boldsymbol\xi f:\mathbb{Z}_{n}\rightarrow \mathbb{Z}}
	\end{equation*}		
	is not well-defined. Therefore $\boldsymbol\xi$ is not a functor of relations.
\end{example}

If $\textbf{Can}$ and $\overline{\textbf{Can}}$, are categories of canonical hypergroups, which their morphisms are good and inclusion homomorphisms, respectively, then we can prove the following theorem using the results obtained in the previous section.

\begin{theorem}
	Every hom-invariant sequence $\boldsymbol\theta$ of regular relations, is a functor on $\overline{\textbf{Can}}$ (or $\textbf{Can}$).
	\begin{proof}
		Let $f:H\rightarrow K$ be an inclusion homomorphism of canonical hypergroups. Then by Theorem \ref{2.2}, $\boldsymbol\theta(H)=H/\theta_{H}$ and $\boldsymbol\theta(K)=H/\theta_{K}$ are canonical hypergroups. Also $\boldsymbol\theta f:\boldsymbol\theta(H)\rightarrow \boldsymbol\theta(H)$ by $\boldsymbol\theta f(\theta_{H}(h))=\theta_{K}(f(h))$ is well defined, and 
		\begin{align*}
			\boldsymbol\theta f(\theta_{H}(h_{1})+\theta_{H}(h_{2}))&=\boldsymbol\theta f(\theta_{H}(h_{1}+h_{2}))= \theta_{K}(f(h_{1}+h_{2}))\\
			&\subseteq \theta_{K}(f(h_{1})+f(h_{2}))=\theta_{K}(f(h_{1}))+\theta_{K}(f(h_{2}))\\
			&= \boldsymbol\theta f(\theta_{H}(h_{1}))+\boldsymbol\theta f(\theta_{H}(h_{2})).
		\end{align*}
		Hence $\boldsymbol\theta\in Hom(\boldsymbol\theta(H),\boldsymbol\theta(K))$. Moreover, if $g\in Hom_{_{\overline{\textbf{Can}}}}(K,T)$, then $g\circ f$ is hom-invariant, because $g\circ f(S_{H})\subseteq S_{T}$. Therefore
		\begin{center}
			$\boldsymbol\theta(g\circ f)(\theta_{H}(h))=\theta_{T}(g(f(h)))=\boldsymbol\theta g(\theta_{K}(f(h)))= \boldsymbol\theta g(\boldsymbol\theta f(\theta_{H}(h)))=\boldsymbol\theta g\circ \boldsymbol\theta f(\theta_{H}(h))$.
		\end{center}
	\end{proof}	
\end{theorem}

Moreover, if $\textbf{F.S.R}(\overline{\textbf{H}})$ is the set of all hom-invariant sequences of strongly regular relations on regular hypergroups, then the following operations, give it a lattice structure, \cite[Lemma 4.1]{Afshar-Ameri 2}.
\begin{center}			
	$\big\{ \theta_{_{H}}\big\}_{H\in I}\vee\big\{\delta_{_{H}}\big\}_{H\in I}=\big\{ \theta_{_{H}}\vee\delta_{_{H}}\big\}_{H\in I} \ \ \ , \ \ \ \big\{ \theta_{_{H}}\big\}_{H\in I}\cap\big\{\delta_{_{H}}\big\}_{H\in I}=\big\{ \theta_{_{H}}\cap\delta_{_{H}}\big\}_{H\in I}$.						
\end{center}				
Because, for every homomorphism $f:H\rightarrow K$ of regular hypergroups, we have $f(S_{\theta_{H}})\subseteq S_{\theta_{K}}$ and $f(S_{\delta_{H}})\subseteq S_{\delta_{K}}$. Hence by \cite[Proposition 4.3]{Afshar-Ameri 2}:
\begin{center}			
	$f(S_{\theta_{H}\vee\delta_{H}})=f(S_{\theta_{H}}+S_{\delta_{H}})\subseteq f(S_{\theta_{H}})+f(S_{\delta_{H}})\subseteq S_{\theta_{K}}+S_{\delta_{K}}=S_{\theta_{K}\vee\delta_{K}}$.					
\end{center} 
We call the members of $\textbf{F.S.R}(\overline{\textbf{H}})$, \emph{functor of relations}. Let $\textbf{F.R}(\overline{\textbf{Can}})$ be the set of all hom-invariant sequences of regular relations on canonical hypergroups. Then by Proposition \ref{3.3} and Theorem \ref{T3.6}, $\textbf{F.S.R}(\overline{\textbf{Can}})$ is sublattice (full subcategory) of $\textbf{F.R}(\overline{\textbf{Can}})$.

\begin{definition} \label{47}
	A functor of relations $\boldsymbol\theta:\mathcal{C}\rightarrow \mathcal{D}$ is said to be \emph{full surjective} if it is full and $\boldsymbol\theta(ob(\mathcal{C}))=ob(\mathcal{D})$. Also $H\in ob(\mathcal{C})$ is said to be $\boldsymbol\theta$-\emph{universal}, if $\boldsymbol\theta(H)$ is universal object in $\mathcal{D}$.
\end{definition}

\begin{example}
	Clearly, $\boldsymbol\beta=\{\beta_{H}\}_{H\in I}$, $\boldsymbol\gamma=\{\gamma_{H}\}_{H\in I}$, $\boldsymbol\eta=\{\eta_{a,H}^{\ast}\}_{H\in I}$, $\boldsymbol\tau=\{\tau_{n,H}^{\ast}\}_{H\in I}$, $\boldsymbol\Delta=\{\Delta_{H}\}_{H\in I}$ and $\boldsymbol\nabla=\{\nabla_{H}\}_{H\in I}$ are functors of relations from $\overline{\textbf{H}}$ to $\textbf{Group}$. Also, $\boldsymbol\Delta$ and $\boldsymbol\nabla$ are full surjective, but the others are not. 
	
\end{example}

\begin{example}
	Let $\overline{\textbf{H.G}}$ be a category which its morphisms are inclusion homomorphisms and its objects are abelian regular hypergroups such as $H$, such that $E_{H}\subseteq x+z$ and $C(x+y)=C(x)+C(y)$ for every $x,y\in H$ and $z\in C(x)$. Now cosider the following relation on hypergroup $H$:
	\begin{equation*}
		(x,y)\in \lambda_{H} \Longleftrightarrow C(x)=C(y).
	\end{equation*}	
	If $C(x)=C(y)$, then $C(x+z)=C(x)+C(z)=C(y)+C(z)=C(y+z)$, for every $z\in H$. Hence for every $s\in x+z$, there is a $t\in y+z$ such that $C(s)=C(t)$. So $\lambda_{H}$ is regular. Also $\lambda_{H}(x)=C(x)$. Let $H$ be the hypergroup defined in Example \ref{3--7}, then $H/\lambda_{H}$ is the following hypergroup:
	\begin{center}
		\begin{tabular}{c|cccccc}		
			$+$& $\{0\}$  & $\{\pm x\}$  &$\{y\}$  &$\{\pm z\}$  &$\{u\}$  &$\{v\}$ \\
			\hline
			$\{0\}$& $\{0\}$  & $\{\pm x\}$        &$\{y\}$       &$\{\pm z\}$      &$\{u\}$      &$\{v\}$ \\
			$\{\pm x\}$& $\{\pm x\}$&$\{0\},\{u\},\{v\}$  & $\{\pm x\},\{\pm z\}$  &$\{u\},\{v\}$   &$\{\pm x\},\{\pm z\}$   &$\{\pm x\}$  \\
			$\{y\}$& $\{y\}$  & $\{\pm x\},\{\pm z\}$    &$\{0\},\{u\}$   &$\{\pm x\}$      &$\{y\},\{v\}$  &$\{u\}$   \\
			$\{\pm z\}$& $\{\pm z\}$  & $\{y\},\{u\}$   &$\{\pm x\}$    &$\{0\},\{v\}$   &$\{\pm x\}$   &$\{\pm z\}$ \\
			$\{u\}$& $\{u\}$  & $\{\pm x\},\{\pm z\}$    &$\{y\},\{v\}$   &$\{\pm x\}$      &$\{0\},\{u\}$  &$\{y\}$    \\
			$\{v\}$& $\{v\}$  & $\{\pm x\}$        &$\{u\}$       &$\{\pm z\}$      &$\{y\}$      &$\{0\}$    
		\end{tabular}
	\end{center}	
	Also, $\lambda_{\mathfrak{R}}=\beta_{\mathfrak{R}}$ and $\mathfrak{R}/\lambda_{\mathfrak{R}}\cong \mathbb{Z}$. If $\boldsymbol{\lambda}=\{\lambda_{H}\}_{H\in I}$, $f\in Hom(H,K)$ and $(x,y)\in \lambda_{H}$, then $C(x)=C(y)$ and $f(C(x))=f(C(y))$. Clearly, $f(C(x))\subseteq C(f(x))$, and $f(C(y))\subseteq C(f(y))$, therefore $C(f(x))\cap C(f(y))\neq\emptyset$, and $C(f(x))=C(f(y))$. So, we have $(f(x),f(y))\in \lambda_{K}$, which means $f(S_{H})\subseteq S_{K}$. Therefore, $\boldsymbol{\lambda}=\{\lambda_{H}\}_{H\in I}$ is a functor from the category $\overline{\textbf{H.G}}$ to itself.
	
\end{example}

\section{Applications}\label{Sec5}

From now on, by hypergroup, we mean the abelian regular hypergroups such as $H$, which  $E_{H}\subseteq x+z$ and $C(x+y)=C(x)+C(y)$ for every $x,y\in H$ and $z\in C(x)$. We use $\textbf{H.G}$ to denote the category of these hypergroups with good homomorphisms, and if the hypergroups are strongly regular as well, we use $\textbf{S.H.G}$.\\
The bar symbol on the category of hypergroups, indicates that the morphisms are inclusion homomorphisms, that map identity elements to identity elements, i.e., $f(h_{1}+h_{2})\subseteq f(h_{1})+f(h_{2})$ and $f(E_{H})\subseteq E_{K}$, for every $f\in Hom(H,K)$ and $h_{1},h_{2}\in H$.\\
So,		
$\textbf{Can}<\textbf{S.H.G}<\textbf{H.G}$, and 		
$\overline{\textbf{Can}}<\overline{\textbf{S.H.G}}<\overline{\textbf{H.G}}$. Where $<$ means full subcategory.

\begin{definition}
	If $\boldsymbol\theta\in \textbf{F.S.R}(\overline{\textbf{H.G}})$, then the hypergroup $H$ is called the $\boldsymbol\theta$-free hypergroup, whenever $\boldsymbol\theta(H)$ is free in category of abelian groups $\textbf{A.Group}$.
\end{definition}

\begin{example} \label{ex}
	Assume that $\mathfrak{R}$ is the hypergroup introduced in Example \ref{EX5}, and $H$ is the hypergroup introduced in Example \ref{3--7}. Since $\boldsymbol\beta(\mathfrak{R})=\mathbb{Z}$, then $\mathfrak{R}$ is $\boldsymbol\beta$-free because $\mathbb{Z}$ is free in $\textbf{A.Group}$. Also, $\mathfrak{R}\times\mathfrak{R}$, $\mathfrak{R}\times\mathbb{Z}$,... are $\boldsymbol\beta$-free.
	But $\boldsymbol\beta(H)\cong \mathbb{Z}_{4}$, and $H$ is not $\boldsymbol\beta$-free, because $\mathbb{Z}_{4}$ is not free in $\textbf{A.Group}$.
\end{example}

Let $(R,\oplus,\odot)$ be a hyperring, where $0_{R}\neq\emptyset$. A hypergroup $(H,+)$ together with the scalar multiplication $\cdot:R\times H\longrightarrow\mathcal{P}^{\ast}(H)$ is called a \emph{left $R$-hypermodule} if for all $r,s\in R$ and $x,y\in H$ the following axioms holds:
\begin{itemize}
	\item[($i$)] $r\cdot(x+y)=r\cdot x+r\cdot y$ ;
	\item[($ii$)] $(r+s)\cdot x\subseteq r\cdot x+s\cdot y$ ;
	\item[($iii$)] $(rs)\cdot x\subseteq r\cdot(s\cdot x)$ ;
	\item[($iv$)] $0_{R}\cdot x=E_{H}$.
\end{itemize}
If $R$ is a Krasner hyperring, and $H$ is a canonical hypergroup, then $H$ is called a \emph{Krasner $R$-hypermodule}. For strongly regular hypergroups, in ($iii$), equality will be established.\\
To refer to the categories of Krasner $\mathbb{Z}$-hypermodules and $R$-hypermodules, we use the symbols $\mathbb{Z}\textbf{-K.H.M}$ and $R\textbf{-H.M}$, respectively.

\begin{lemma}	\label{44}
	Every hypergroup is $\mathbb{Z}$-hypermodule.
	\begin{proof}
		Consider the following scalar multiplication:
		\begin{equation*}\label{eq5}
			\cdot:\mathbb{Z}\times H\rightarrow \mathcal{P}^{\ast}(H)	
		\end{equation*}
		\begin{equation}
			n\cdot h=\begin{cases}
				nh&    n>0 \\
				E_{H}&          n=0 \\
				C(-nh)&         n<0
				
			\end{cases}		
		\end{equation}
		Where $nh=\sum_{i=1}^{n}h$ and $C(-nh)=\sum_{i=1}^{-n}C(h)$. Let $m,n\in\mathbb{Z}$ and $x,y\in H$, then:\\
		(i) If $n>0$, then $n\cdot(x+y)=n\cdot x+n\cdot y$ and for $n<0$ we have $n\cdot(x+y)=C(-n(x+y))=C(-nx+(-n)y)=C(-nx)+C(-ny)=n\cdot x+n\cdot y$. If $n=0$, then $0\cdot(x+y)=E_{H}$.\\
		(ii) If $m,n>0$, then $(m+n)\cdot x=m\cdot x+n\cdot x$. If $m,n>0$ but $n<0$, then:
		\begin{center}
			$m\cdot x+n\cdot x=mx+C(-nx)=mx+\sum_{i=1}^{-n}C(x)=(m+n)x+(-n)x+\sum_{i=1}^{-n}C(x).$
		\end{center}
		Since $E_{H}\subseteq(-n)x+\sum_{i=1}^{-n}C(x)$, then $(m+n)\cdot x\subseteq m\cdot x+n\cdot x$. If $m+n=0$, then $n=-m$ and $(m-m)x=0\cdot x=E_{H}\subseteq mx-mx$. Let $m+n<0$, then $(m+n)\cdot x=C(-(m+n)x)=C(-mx)+C(-nx)=m\cdot x+n\cdot x.$. If $m+n<0$ and $n>0$, then $(m+n)\cdot x=C(-(m+n)x)$ and
		\begin{center}
			$m\cdot x+n\cdot x=C(-mx)+nx=\sum_{i=1}^{-m}C(x)+nx=(m+n)x+(-n)x+nx.$
		\end{center}
		Since $E_{H}\subseteq-nx+nx$, then $(m+n)\cdot x\subseteq m\cdot x+n\cdot x$.\\
		(iii) If $m+n>0$, then $(mn)\cdot x=m\cdot(n\cdot x)$. Let $m>0$ and $n<0$, then $(mn)\cdot x=C(-(mn)x)$ and $m\cdot(n\cdot x)=m\cdot(nx)=C(-m(nx))=C(-mn\cdot x)=mn\cdot x$. For $m,n<0$ we have $(mn)\cdot x=mnx$ and $m\cdot (n\cdot x)=m\cdot C(-nx)=C(-mc(-nx))=C(C(mnx))$. So $mnx\subseteq C(C(mnx))$. If $mn=0$, let $n=0$. Then $(mn)\cdot x=0\cdot x=E_{H}$ and $m\cdot(n\cdot x)=m\cdot E_{H}=\sum_{i=1}^{m}E_{H}$. Since $E_{H}\subseteq\sum_{i=1}^{m}E_{H}$ then $(mn)\cdot x\subseteq m\cdot(n\cdot x)$.
	\end{proof}	
\end{lemma}

\begin{example}
	Consider homomorphism $f\in Hom(\mathfrak{R},\mathbb{Z})$ by $(r,s)\mapsto 2r$, in category $\textbf{H.G}$. Then $\boldsymbol\beta f\in Hom(\mathbb{Z},\mathbb{Z})$, where $\boldsymbol\beta f(r)=2r$. 
\end{example}

For a fixed hyperring $R$, we use the bar symbol on the category of $R$-hypermodules to indicate that the morphisms are inclusion homomorphisms, and if $f: M\rightarrow N$ is an inclusion homomorphism of $R$-hypermodules, then $f(r\cdot m)\subseteq r\cdot f(m)$, for every $r\in R$ and $m\in M$.

In Example \ref{ex}, if $(a,b)\odot(c,d)=\big\{(ac,x);\ x\in\mathbb{Z} \big\}$, then $(\mathfrak{R},\oplus,\odot)$ is a general hyperring and the following Lemma is obtained.

\begin{lemma}	\label{45}
	Every hypergroup is $\mathfrak{R}$-hypermodule.
	\begin{proof}	
		Let $H$ be a hypergroup. It is enough to consider the scalar multiplication $\circ:\mathfrak{R}\times H\rightarrow \mathcal{P}^{\ast}(H)$ as $(r,s)\circ h=r\cdot h$, for every $h\in H$ and $(r,s)\in\mathfrak{R}$. Given that $H$ is $\mathbb{Z}$-hypermodule, the rest of the proof is straightforward.	
	\end{proof}	
\end{lemma}

\begin{proposition}
	\emph{$\textbf{Can}=\mathbb{Z}\textbf{-K.H.M}$} and 
	\emph{$\overline{\textbf{Can}}=\mathbb{Z}\textbf{-}\overline{\textbf{K.H.M}}$}.
	\begin{proof}	
		Let $f\in Hom(H,K)$ be a good homomorphism of canonical hypergroups. Then for every $n\in \mathbb{N}$ and $h\in H$; $f(n\cdot h)=f(h+...+h)=n\cdot f(h)$, and $f(0_{H})=0_{K}$. Since $f(0_{H})\in f(-h+h)=f(-h)+f(h)$, then $f(-h)=-f(h)$. So $f(h-h)=f(h+(-h))=f(h)+f(-h)=f(h)-f(h)$. Therefore $f(n\cdot h)=n\cdot f(h)$ for every $n\in \mathbb{Z}$. Hence $Hom(H,K)=Hom_{\mathbb{Z}}(H,K)$.
		This proof can be similarly stated for inclusion homomorphisms.
	\end{proof}	
\end{proposition}

\begin{theorem}
	\emph{$\overline{\textbf{H.G}}=\mathbb{Z}\textbf{-}\overline{\textbf{H.M}}=\mathfrak{R}\textbf{-}\overline{\textbf{H.M}}$}.
	\begin{proof}
		Let $f\in Hom(H,K)$ be an inclusion homomorphism of hypergroups. and $n\in \mathbb{Z}$. If $n>0$, then $f(n\cdot h)=f(h+...+h)\subseteq f(h)+...+f(h)=n\cdot f(h)$. If $n=0$, then $f(0\cdot h)=f(E_{H})\subseteq E_{K}=0\cdot f(h)$.\\
		Let $z\in C(h)$, then $E_{H}\subseteq z+h$ and $f(E_{H})\subseteq f(z)+f(h)$. So $E_{K}\cap f(z)+f(h)\neq \emptyset$	and $f(z)\in C(f(h))$. Therefore $f(C(h))\subseteq C(f(h))$. If $n<0$, then:
		\begin{align*}
			f(n\cdot h)&=f(C(-nh))=f(\sum_{i=1}^{-n}C(h))\\
			&=\sum_{i=1}^{-n}f(C(h))\subseteq \sum_{i=1}^{-n}C(f(h))\\
			&=C(\sum_{i=1}^{-n}f(h))=C(-nf(h))\\
			&=n\cdot f(h).
		\end{align*}
		Conversely, if $n=0$, then $f(E_{H})=f(0\cdot h)\subseteq 0\cdot f(h)=E_{K}$. Therefore $Hom(H,K)=Hom_{\mathbb{Z}}(H,K)$. Now, by Lemma \ref{44} we have $\overline{\textbf{H.G}}=\mathbb{Z}\textbf{-}\overline{\textbf{H.M}}$.\\
		Using Lemma \ref{45} and since $(r,s)\circ h=r\cdot h$, for every $(r,s)\in\mathfrak{R}$, the expression $\overline{\textbf{H.G}}=\mathfrak{R}\textbf{-}\overline{\textbf{H.M}}$ is also proved in a similar way.
	\end{proof}		
\end{theorem}

Consider the $\mathbb{Z}$-hypermodules $H$, $K$ and hypergroup $C$. A function $f:H\times K\rightarrow C$ is said to be an \emph{inclusion middle linear map}, whenever:
\begin{itemize}
	\item[($i$)] $f(h_{1}+h_{2},k)\subseteq f(h_{1},k)+f(h_{2},k)$;
	\item[($ii$)] $f(h,k_{1}+k_{2})\subseteq f(h,k_{1})+f(h,k_{2})$;
	\item[($iii$)] $f(nh,k)=f(h,nk)$.
\end{itemize}
For every $h,h_{1},h_{2}\in H$, $k,k_{1},k_{2}\in K$ and $n\in \mathbb{Z}$. Moreover, if condition $(iii)$ is such that $f(nh,k)=f(h,nk)\subseteq nf(h,k)$, then $f$ is called \emph{inclusion bilinear map}.\\
By $\bar{\mathcal{M}}(H,K)$, we mean the category of inclusion middle linear maps on $H\times K$, that for every inclusion middle linear maps $f_{i}:H\times K\rightarrow C_{i}$ and $f_{j}:H\times K\rightarrow C_{j}$; $Hom(f_{i},f_{j})=\{g\in Hom_{\overline{\textbf{H.G}}}(C_{i},C_{j});\ g\circ f_{i}=f_{j}\}$. Indeed, morphisms are all homomorphisms of hypergroups in $\overline{\textbf{H.G}}$, such as $g$, such that the following diagram is commutative:
\begin{center}
	\begin{tikzcd}	
		H\times K  \ar[d,"f_{j}"'] \ar[r,"f_{i}"]& C_{i} \ar[dl,"g"] \\
		C_{j}
	\end{tikzcd}
\end{center}
From now on, we will work with functors like $\boldsymbol\theta$ that have the property $\boldsymbol\theta(\prod_{j\in J} H_{j})=\prod_{j\in J}\boldsymbol\theta(H_{j})$.

\begin{definition}
	A sequence of relations $\boldsymbol\theta=\big\{ \theta_{_{H}}\big\}_{H\in I}$ is called \emph{linear-invariant} if $f(S_{H\times K})\subseteq S_{C}$, for every hypergroups $H$, $K$, $C$ and every inclusion middle linear map $f:H\times K\rightarrow C$. Also, $\boldsymbol\theta$ is said to be \emph{invariant} if it is hom-invariant and linear-invariant.
	
\end{definition}

\begin{theorem}
	Every invariant	sequence of relations, is a functor from $\bar{\mathcal{M}}(H,K)$ to the category of middle linear maps $\mathcal{M}(\boldsymbol\theta(H),\boldsymbol\theta(K))$.
	\begin{proof}
		Consider the inclusion middle linear map $f:H\times K\mapsto C$ and the invariant sequence of relations $\boldsymbol\theta$. We must show that
		\begin{center}
			$\underset{(\boldsymbol\theta(h),\boldsymbol\theta(k)) \mapsto \ \boldsymbol\theta(f(h,k)) }
			{\boldsymbol\theta f:\boldsymbol\theta(H)\times\boldsymbol\theta(K) \rightarrow \boldsymbol\theta(C)}$
		\end{center}
		is a middle linear map of abelian groups.\\
		(i) If $(\boldsymbol\theta(h_{1}),\boldsymbol\theta(k_{1}))=(\boldsymbol\theta(h_{2}),\boldsymbol\theta(k_{2}))$, then $(h_{1},k_{1})+S_{H\times K}=(h_{2},k_{2})+S_{H\times K}$. So $f((h_{1},k_{1})+S_{H\times K})=f((h_{2},k_{2})+S_{H\times K})$. Therefore
		\begin{align*}
			f((h_{2},k_{2})+S_{H\times K})&\subseteq f(h_{2},k_{2})+f(S_{H\times K})\subseteq f(h_{2},k_{2})+S_{C}\\
			&=\boldsymbol\theta(f(h_{2},k_{2}))=\boldsymbol\theta f(\boldsymbol\theta(h_{2}),\boldsymbol\theta(k_{2})).
		\end{align*}
		Hence $\boldsymbol\theta f(\boldsymbol\theta(h_{1}),\boldsymbol\theta(k_{1}))=\boldsymbol\theta f(\boldsymbol\theta(h_{2}),\boldsymbol\theta(k_{2}))$, and $\boldsymbol\theta f$ is well defined.\\
		(ii) Since $f(h_{1}+h_{2},k)\subseteq f(h_{1},k)+f(h_{2},k)$ and $\boldsymbol\theta$ is a sequence of strongly regular relations, then 
		\begin{center}
			$\boldsymbol\theta\big(f(h_{1}+h_{2},k)\big)=\boldsymbol\theta\big(f(h_{1},k)+f(h_{2},k)\big)=\boldsymbol\theta\big(f(h_{1},k)\big)+\boldsymbol\theta\big(f(h_{2},k)\big)$.	
		\end{center}
		Hence
		\begin{center}
			$\boldsymbol\theta f\big(\boldsymbol\theta(h_{1})+\boldsymbol\theta(h_{2}),\boldsymbol\theta(k)\big)
			=\boldsymbol\theta f\big(\boldsymbol\theta(h_{1}),\boldsymbol\theta(k)\big)+\boldsymbol\theta f\big(\boldsymbol\theta(h_{2}),\boldsymbol\theta(k)\big)$.		
		\end{center}
		Similarly, the following statements are also proven:
		\begin{center}
			$\boldsymbol\theta f\big(\boldsymbol\theta(h),\boldsymbol\theta(k_{1})+\boldsymbol\theta(k_{2})\big)
			=\boldsymbol\theta f\big(\boldsymbol\theta(h),\boldsymbol\theta(k_{1})\big)+\boldsymbol\theta f\big(\boldsymbol\theta(h),\boldsymbol\theta(k_{2})\big)$;
		\end{center}	
		\begin{center}	
			$\boldsymbol\theta f\big(\boldsymbol\theta(h),n\boldsymbol\theta(k)\big)
			=\boldsymbol\theta f\big(n\boldsymbol\theta(h),\boldsymbol\theta(k)\big)$.			
		\end{center}
		(iii) If $f_{i}:H\times K\mapsto C_{i}$ and $f_{j}:H\times K\mapsto C_{j}$ are inclusion middle linear maps and $g\in Hom(f_{i},f_{j})$, then $g\circ f_{i}=f_{j}$. Since $\boldsymbol\theta$ is invariant, then $g(S_{C_{i}})\subseteq S_{C_{j}}$, and by Theorem \ref{T.1}. $\boldsymbol\theta g:\boldsymbol\theta(C_{i})\rightarrow\boldsymbol\theta(C_{i})$ is a homomorphism of abelian groups. On the other hand, according to part (ii), $\boldsymbol\theta f_{i}$ and $\boldsymbol\theta f_{j}$ are middle linear maps. Also. by Theorem \ref{T.1}, $\boldsymbol\theta$ is a functor from $\overline{\textbf{H.G}}$ to $\textbf{A.Group}$. Hence $\boldsymbol\theta(g\circ f_{i})=\boldsymbol\theta g\circ\boldsymbol\theta f_{i}$. But $\boldsymbol\theta(g\circ f_{i})=\boldsymbol\theta f_{j}$, and therefore the following diagram is commutative:
		\begin{center}
			\begin{tikzcd}	
				\boldsymbol\theta(H)\times\boldsymbol\theta(K) \ar[r,"f_{i}"] \ar[rd,"f_{j}"']   & \boldsymbol\theta(C_{i}) \ar[d,"\boldsymbol\theta g"] \\
				& \boldsymbol\theta(C_{j})
			\end{tikzcd}
		\end{center}
		(iv) If $g_{1}\in Hom(f_{i},f_{j})$ and $g_{2}\in Hom(f_{j},f_{k})$ then $g_{1}\circ f_{i}=f_{j}$ and $g_{2}\circ f_{j}=f_{k}$. Thus $g_{2}\circ g_{1}\circ f_{i}=f_{k}$. So $g_{2}\circ g_{1}\in Hom(f_{i},f_{k})$ and $\boldsymbol\theta(g_{2}\circ g_{1})\in Hom(\boldsymbol\theta f_{i},\boldsymbol\theta f_{k})$. Similarly $\boldsymbol\theta(g_{2})\circ \boldsymbol\theta(g_{1})\in Hom(\boldsymbol\theta f_{i},\boldsymbol\theta f_{k})$. Now by Theorem \ref{T.1}, $\boldsymbol\theta(g_{2}\circ g_{1})=\boldsymbol\theta(g_{2})\circ \boldsymbol\theta(g_{1})$.\\
		Therefore $\boldsymbol\theta$ is a functor from $\bar{\mathcal{M}}(H,K)$ to  $\mathcal{M}(\boldsymbol\theta(H),\boldsymbol\theta(K))$.
	\end{proof}
\end{theorem}

\begin{remark}
	We know that the universal object of the category $\mathcal{M}\big(\boldsymbol\theta(H),\boldsymbol\theta(K)\big)$, is $\otimes:\boldsymbol\theta(H)\times\boldsymbol\theta(K)\rightarrow\boldsymbol\theta(H)\otimes_{\mathbb{Z}}\boldsymbol\theta(K)$, \cite{Hungerford}. On the other hand, by Definition \ref{47}, the $\boldsymbol\theta$-universal objects of $\bar{\mathcal{M}}(H,K)$ are all inclusion middle linear maps such as $\boldsymbol\otimes:H\times K\rightarrow C$, that $\boldsymbol\theta\boldsymbol\otimes=\otimes$. Therefore: 
	\begin{equation*}\label{eq5}
		H\boldsymbol\otimes_{_{\mathbb{Z}}}K=\{C ;\ \boldsymbol\theta(C)\cong\boldsymbol\theta(H)\otimes_{_{\mathbb{Z}}}\boldsymbol\theta(K)\}
	\end{equation*}
	which we call it the $\boldsymbol\theta$-\emph{tensor product} of hypergroups $H$ and $K$.
	Moreover $\boldsymbol\theta(H\boldsymbol\otimes_{_{\mathbb{Z}}}K)\cong\boldsymbol\theta(H)\otimes_{_{\mathbb{Z}}}\boldsymbol\theta(K)$.
	
\end{remark}

\begin{example}			
	Consider the hypergroups $H$ and $\mathfrak{R}$, as defined in Example \ref{3--7} and \ref{EX5}, respectively.\\
	(i) $\mathbb{Q} \boldsymbol\otimes_{_{\mathbb{Z}}}^{\boldsymbol\beta}\mathbb{Z}_{n}=\{T; \boldsymbol\beta(T)\cong\mathbb{Q} \otimes_{_{\mathbb{Z}}}\mathbb{Z}_{n}\}=\{T;\boldsymbol\beta(T)\cong \mathbb{Z}_{1}\}=all \ \boldsymbol\beta-trivial\ hypergroups$.\\
	(ii) $H \boldsymbol\otimes_{_{\mathbb{Z}}}^{\boldsymbol\beta}\mathbb{Z}_{6}=\{T; \boldsymbol\beta(T)\cong\boldsymbol\beta(H) \otimes_{_{\mathbb{Z}}}\mathbb{Z}_{6}\}=\{T; \boldsymbol\beta(T)\cong\mathbb{Z}_{4} \otimes_{_{\mathbb{Z}}}\mathbb{Z}_{6}\}=\{T;\boldsymbol\beta(T)\cong \mathbb{Z}_{2}\}$.\\
	(iii) $\mathfrak{R} \boldsymbol\otimes_{_{\mathbb{Z}}}^{\boldsymbol\beta}\mathbb{Z}_{2}=\{T; \boldsymbol\beta(T)\cong\boldsymbol\beta(\mathfrak{R}) \otimes_{_{\mathbb{Z}}}\mathbb{Z}_{2}\}=\{T; \boldsymbol\beta(T)\cong\mathbb{Z} \otimes_{_{\mathbb{Z}}}\mathbb{Z}_{2}\}=\{T;\boldsymbol\beta(T)\cong \mathbb{Z}_{2}\}$.
	
\end{example}

\section{Conclusion}\label{sec6}

Therefore, we can consider every regular relation as a congruence relation. Moreover, strongly regular relations are actually a specific type of these congruence relations, which are modulo subhypergroups containing $S_{\beta}$. One of our future investigations will be to generalize the results obtained here to any arbitrary hypergroup.\\
We also showed that every sequence of strongly regular relations, is functorial. Therefore, by considering a specific property, we can construct a functor using a sequence of strongly regular relations that preserves that property.\\
For example, if $\theta^{\ast}$ is the smallest strongly regular relation such that the quotients obtained from it are cyclic group, it is sufficient to choose the sequence $\boldsymbol\theta=\{\theta_{H}\}_{H\in I}$ such that $\theta^{\ast}\subseteq \theta_{H}$ for all $H\in I$.



\end{document}